\documentclass{article}

\usepackage{authblk}

\usepackage{geometry}
\usepackage{setspace}
\usepackage{amsmath}
\usepackage{amssymb}
\usepackage{mathtools}
\usepackage{bm}
\usepackage{blkarray} %block array for labelling matrix columns

\usepackage{hyperref}
\usepackage{xcolor}
\definecolor{darkgreen}{RGB}{0,127,173}
\definecolor{darkblue}{RGB}{0,0,200}
\hypersetup{
	colorlinks,
        linkcolor={darkblue},
        citecolor={darkgreen},
        urlcolor={darkblue}
}
\usepackage{verbatim}

\usepackage{longtable}
\allowdisplaybreaks

\usepackage{graphicx}
\usepackage{subfig}
\usepackage{natbib}
\title{A compact reformulation of the two-stage robust resource-constrained project scheduling problem}
\author[1]{Matthew Bold}
\author[2]{Marc Goerigk}
\affil[1]{STOR-i Centre for Doctoral Training, Lancaster University, United Kingdom}
\affil[2]{Network and Data Science Management, University of Siegen, Germany}

\date{}

\begin{document}

\maketitle

\begin{abstract} 
This paper considers the resource-constrained project scheduling problem with uncertain activity durations. We assume that activity durations lie in a budgeted uncertainty set, and follow a robust two-stage approach, where a decision maker must resolve resource conflicts subject to the problem uncertainty, but can determine activity start times after the uncertain activity durations become known. 

We introduce a new reformulation of the second-stage problem, which enables us to derive a compact robust counterpart to the full two-stage adjustable robust optimisation problem. Computational experiments show that this compact robust counterpart can be solved using standard optimisation software significantly faster than the current state-of-the-art algorithm for solving this problem, reaching optimality for almost 50\% more instances on the same benchmark set.
\end{abstract}


\section{Introduction} 

The resource-constrained project scheduling problem (RCPSP) consists of scheduling a set of activities, subject to precedence constraints and limited resource availability, with the objective of minimising the overall project duration, known as the makespan. Given its practical relevance to a number of industries, including construction \citep{kim2013genetic}, manufacturing \citep{gourgand2008employee}, R\&D \citep{vanhoucke2006scheduling}, and personnel scheduling \citep{drezet2008employee}, the RCPSP and many of its variants have been widely studied since a first model was introduced by \cite{pritsker1969multiproject}. The vast majority of this research, however, has examined the RCPSP under the assumption that the model parameters are known deterministically (for a survey of the deterministic RCPSP, see \cite{artigues2008resource}), but clearly, in practice, large projects are subject to non-trivial uncertainties. For instance, poor weather might delay construction times, uncertain delivery times of parts may delay manufacturing activities, and the duration of research activities are inherently uncertain. As a result, in recent years, increasing attention has been given to the uncertain RCPSP, where scheduling decisions must be made whilst activity durations are unknown.  

There exist two main approaches to solving the uncertain RCPSP. The first is to view the problem as a dynamic optimisation problem where scheduling decisions are made each time new information becomes available according to a scheduling policy \citep{igelmund1983algorithmic, igelmund1983preselective, mohring2000linear}. Most recently, \cite{li2015solving} use approximate dynamic programming to find an adaptive closed-loop scheduling policy for the uncertain RCPSP. 

The second approach aims to proactively develop a robust baseline schedule that protects against delays in the activity durations. \cite{zhu2007two} present a two-stage stochastic programming formulation for building baseline schedules for projects with a single resource. \cite{bruni2015stochastic} present a chance-constraint-based heuristic for constructing robust baseline schedules and \cite{lamas2016purely} introduce a procedure for generating robust baseline schedules that does is independent of later reactive scheduling procedures. For a review of both dynamic and proactive project scheduling, see \cite{herroelen2005project}. 

Although frequently referred to as robust, none of the scheduling methods described above make use of robust optimisation as defined in \cite{ben1998robust, ben1999robust, ben2000robust}. Over the last 20 years, robust optimisation has emerged as an effective framework for modelling uncertain optimisation problems. Unlike stochastic programming, robust optimisation does not require probabilistic knowledge of the uncertain data. Instead the robust optimisation approach only assumes that the uncertain data lie somewhere in a given uncertainty set, and then aims to find solutions that are robust for all scenarios that can arise from that uncertainty set.

The applicability of robust optimisation as a method for solving uncertain optimisation problems has increased following the introduction of adjustable robust optimisation \citep{ben2004adjustable, yanikouglu2019survey}. Adjustable robust optimisation extends static robust optimisation into a dynamic setting, where a subset of the decision variables must be determined under uncertainty, whilst other variables can be adjusted following observations of the uncertain data. As well as accurately modelling the decision process undertaken by many real-world decision-makers, adjustable robust optimisation overcomes the over-conservativeness that restricts the applicability of static robust optimisation models. For extensive surveys on robust optimisation, see \cite{ben2009robust, bertsimas2011theory, gorissen2015practical, goerigk2016algorithm}.

Despite the successful application of robust optimisation in many different fields (see \cite{bertsimas2011theory}), so far robust optimisation has been little used in project scheduling. To the best of our knowledge, to date, only three papers have directly applied robust optimisation in the construction of robust baseline project schedules. \cite{artigues2013robust} present an iterative scenario-relaxation algorithm for the uncertain RCPSP with the objective of minimising the worst-case absolute regret \citep{kouvelis1997robust}. \cite{bruni2017adjustable} introduce a two-stage adjustable robust optimisation model with the objective of minimising the worst-case makespan. For the case of budgeted uncertainty, this model is solved using a Benders'-style decomposition approach \citep{benders1962partitioning}. Most recently, \cite{bruni2018computational} present a computational study of solution methods for solving the two-stage adjustable RCPSP. An additional Benders'-style algorithm is compared against a primal decomposition algorithm, as well as the algorithm presented in \cite{bruni2017adjustable}. The primal decomposition algorithm is shown to be the best performing algorithm for solving the two-stage adjustable RCPSP.

This paper presents a new compact reformulation of the two-stage adjustable robust RCPSP with budgeted uncertainty. Computational experiments show that this compact reformulation can be solved using standard optimisation software significantly faster, and for a much greater number of instances than the current best algorithm for solving this problem.

The remainder of this paper is organised as follows: Section \ref{section:problem_description} introduces the two-stage adjustable robust RCPSP in detail, before Section \ref{section:compact_reformulation} derives a compact reformulation of this problem and computational experiments are presented in Section \ref{section:computational_experiments}. Concluding remarks are made in Section \ref{section:conclusion}. 

\section{The two-stage robust RCPSP} \label{section:problem_description}

\begin{figure}[h]
	\centering
	\includegraphics[scale=0.7]{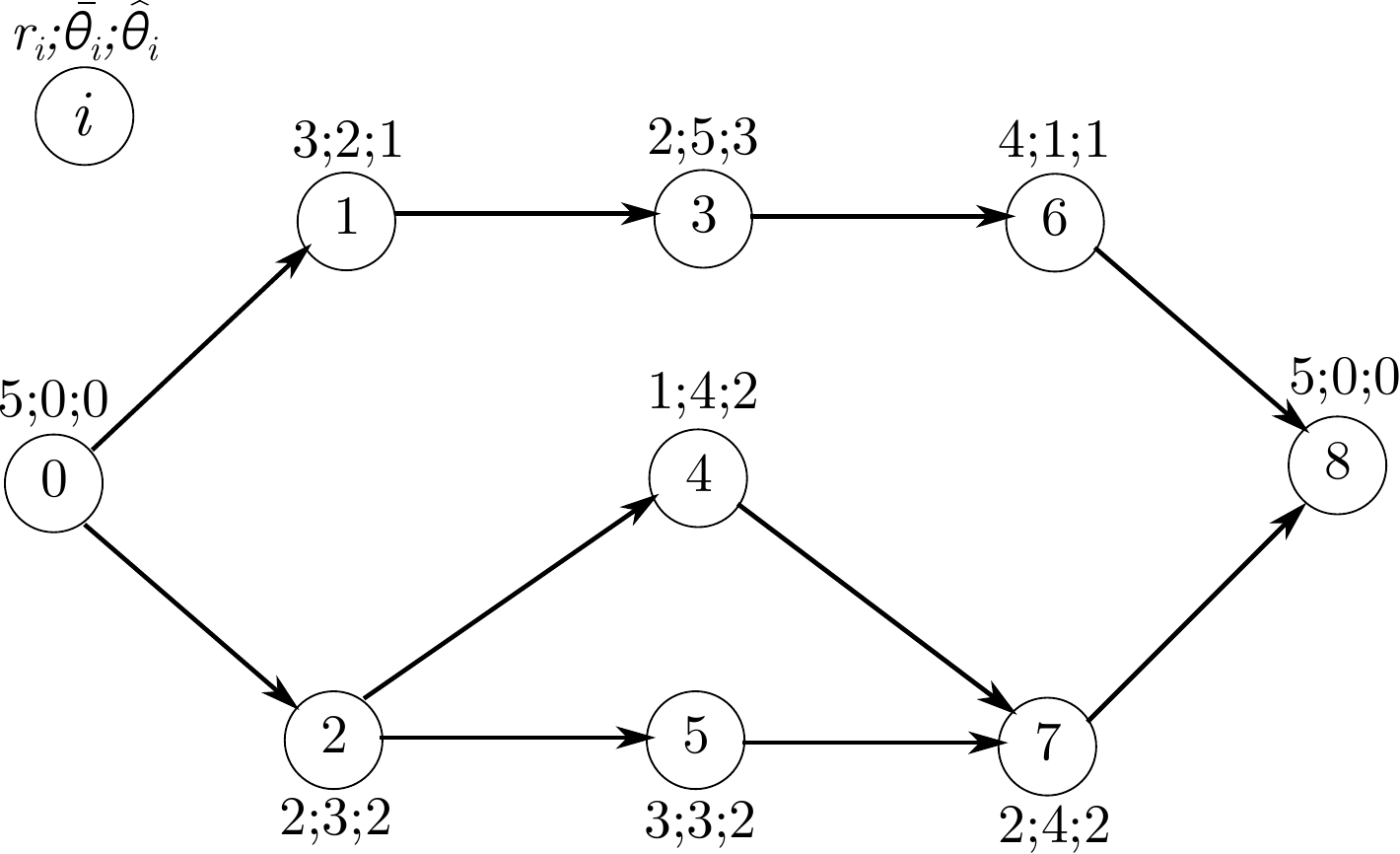}
	\caption{Example project involving seven non-dummy activities and a single resource with $R_1=5$.}
	\label{fig:example_project}
\end{figure}

A project consists of a set $V=\{0,1,\dots,n,n+1\}$ of non-preemptive activities, where $0$ and $n+1$ are dummy source and sink activities with duration 0. Each activity $i\in V$ requires an amount $r_{ik}\geq 0$ of resource $k\in K$, where $K$ is the set of project resource types. Each resource $k\in K$ has a finite availability $R_k$ in each time period. Each activity $i\in V$ has a nominal duration given by $\bar{\theta}_i$, and a worst-case duration given by $\bar{\theta}_i + \hat{\theta}_i$, where $\hat{\theta}_i$ is its maximum deviation. In addition to resource constraints, the project activities must be scheduled in a manner that respects a set $E$ of strict finish-to-start precedence constraints. A project can be represented on a directed graph $G(V,E)$. An example project involving seven non-dummy activities and a single resource is shown in Figure \ref{fig:example_project}. 

We assume that the duration of each activity $i\in V$ lies somewhere between its nominal value $\bar{\theta}_i$ and its worst-case value $\bar{\theta}_i+\hat{\theta}_i$. Additionally, we follow \cite{bertsimas2004price} and assume that only a subset of all activities can simultaneously attain their worst-case values. Hence, the set in which we assume durations can lie, known as the uncertainty set, is given by 
$$\mathcal{U}(\Gamma) = \Bigg\{\theta \in \mathbb{R}^{|V|}_+ \,:\, \theta_i = \bar{\theta}_i + \delta_i\hat{\theta}_i,\,0\leq \delta_i\leq 1\,\forall i\in V,\,\sum_{i\in V}\delta_i\leq \Gamma\Bigg\},$$
 where $\Gamma$ determines the robustness of the solution by controlling the number of activities that are allowed to reach their worst-case duration simultaneously. For $\Gamma = 0$, each activity takes its nominal duration and the problem reduces to the deterministic RCPSP. At the other extreme, when $\Gamma = n$, every activity can take its worst-case duration, and this uncertainty set becomes equivalent to interval uncertainty.

 The robust RCPSP lends itself naturally to a two-stage decision making process, where resource allocation decisions need to be made at the start of the project, before the uncertain activity durations become known, but the activity start times can be decided following the realisation of the activity durations. Hence, resource allocation decisions constitute the set of first-stage decisions, whilst the activity start times constitute the set of second-stage decisions.

\begin{figure}[h]
	\centering
	\includegraphics[scale=0.7]{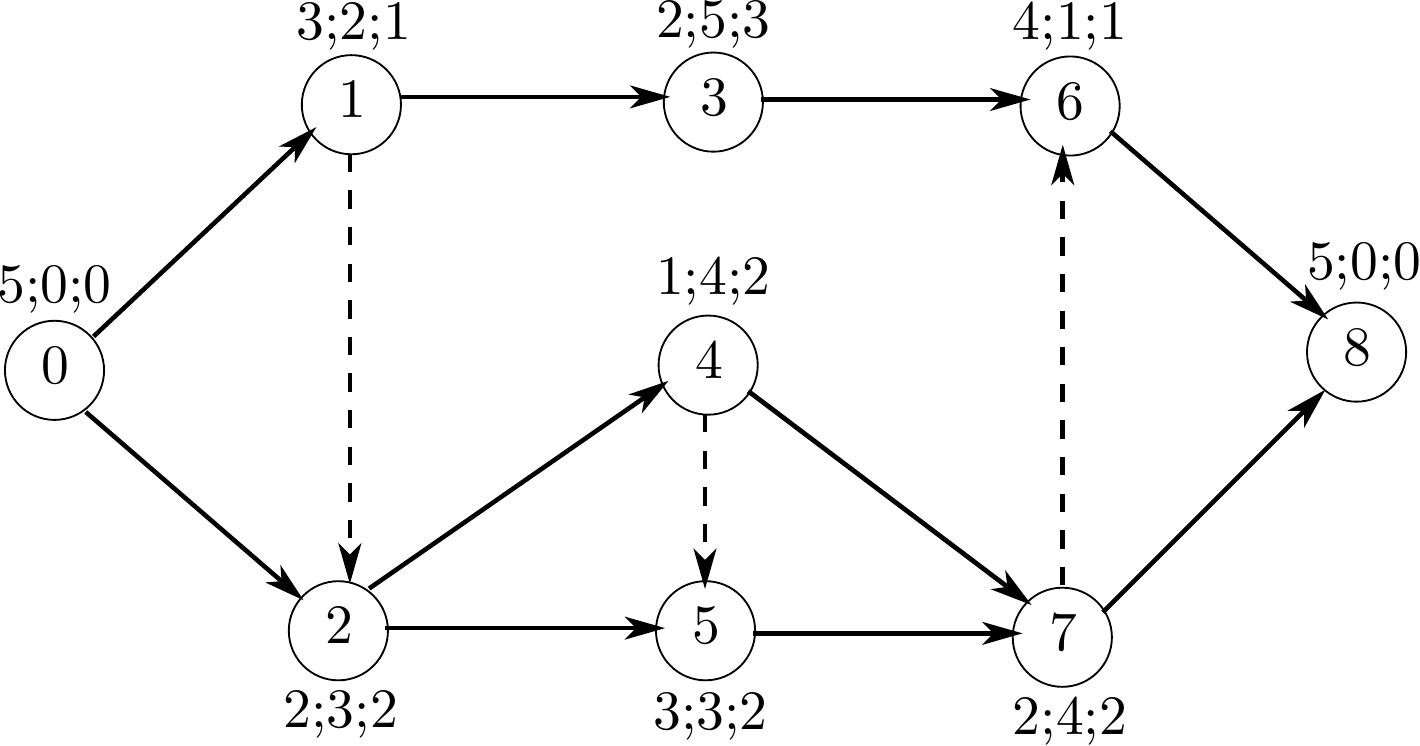}
	\caption{An extension of the example project shown in Figure \ref{fig:example_project}, corresponding to the sufficient selection given by the dashed arcs.}
	\label{fig:sufficient_selection}
\end{figure}

 More specifically, the first-stage resource allocation decisions consist of determining a feasible extension of the project precedence relationships $E$ so that all resource conflicts are resolved. A forbidden set \citep{igelmund1983algorithmic} is any subset $F\subseteq V$ of non-precedence-related activities such that $\sum_{i\in F}r_{ik}>R_k$ for at least one $k\in K$, i.e. the activities of $F$ cannot be processed simultaneously without violating a resource constraint. A minimal forbidden set is a forbidden set that does not contain any other forbidden set as a subset.
 We denote the set of minimal forbidden sets by $\mathcal{F}$. For the example project in Figure \ref{fig:example_project}, $\mathcal{F}=\big\{\{1,5\},\{2,6\},\{5,6\},\{6,7\},\{3,4,5\}\big\}$. The resource conflict represented by each minimal forbidden set can be resolved by adding an additional precedence relationship to the project network. \cite{bartusch1988scheduling} show that solving the RCPSP is equivalent to finding an optimal choice of additional precedence relationships $X\subseteq V^2\setminus E$, such that the extended project network $G'(V,E\cup X)$ is acyclic and contains no forbidden sets. Such an extension $X$ to the project precedence network is referred to as a sufficient selection. Hence, a solution to the first-stage problem corresponds to the choice of a sufficient selection $X$. Figure \ref{fig:sufficient_selection} shows the extended project network for a sufficient selection to the example project shown in Figure \ref{fig:example_project} (arcs in $X$ are dashed).

Given the extended project network resulting from the choice of sufficient selection made in the first stage, the second stage problem consists of determining activity start times in order to minimise the worst-case makespan in this extended network. Since all resource conflicts have been resolved in the first-stage problem, the second stage problem contains no resource constraints.

Hence, the two-stage robust RCPSP under budgeted uncertainty is given by:
\begin{equation}
\min_{X\in \mathcal{X}}\ \max_{\theta\in \mathcal{U}(\Gamma)}\ \min_{S\in\mathcal{S}(X,\theta)} S_{n+1} \label{eqn:problem_eqn3}
\end{equation}
% \begin{align}
% &\min_{X\in \mathcal{X}}\max_{\theta\in \mathcal{U}(\Gamma)}\min_{S(\theta)} S_{n+1}(\theta)\\
% & S_0(\theta) = 0 \quad \forall \theta \in \mathcal{U}(\Gamma)\\
% & S_j(\theta) - S_i(\theta) \geq \theta_i \quad \forall (i,j)\in E\cup X,\,\forall \theta\in\mathcal{U}(\Gamma)\label{eqn:problem_eqn3},
% \end{align}
where $\mathcal{X}$ is the set of sufficient selections and $\mathcal{S}(X,\theta)$ denotes the set of feasible activity start times, depending on activity durations $\theta\in \mathcal{U}(\Gamma)$, as well as on the choice of sufficient selection $X$. We have
\[ \mathcal{S}(X,\theta) = \Bigg\{ S\in\mathbb{R}^{|V|}_+ \,:\, S_0 = 0,\  S_j - S_i \geq \theta_i\,\forall (i,j)\in E\cup X \Bigg\}. \]
% . Note that the second-stage variables $S(\theta)$, representing the activity start times, are dependent on the activity durations $\theta\in \mathcal{U}(\Gamma)$, as well as on the choice of sufficient selection $X$. Their dependence on $X$ is implicit within constraints (\ref{eqn:problem_eqn3}).
To solve this problem we propose a mixed-integer programming formulation, outlined in the following section.

\section{A compact reformulation} \label{section:compact_reformulation}

In this section, we present a compact reformulation of the two-stage robust RCPSP. We begin by first examining the adversarial sub-problem of maximising the worst-case makespan for a given sufficient selection.

\subsection{The adversarial sub-problem}

Suppose the solution to the first-stage problem provides a sufficient selection $X\in \mathcal{X}$, and is given by a vector $y\in\{0,1\}^{V\times V}$
% $\big\{y_{ij} : (i,j)\in V^2\big\}$, 
where 
$$y_{ij}=
\begin{cases}
	1 \quad\textnormal{if } (i,j)\in E\cup X\\
	0 \quad\textnormal{otherwise}.
\end{cases}
$$ 
The second-stage sub-problem that arises can be considered from the point of view of an adversary who wishes to choose the worst-case scenario of delays for the given first-stage solution. Following the adversary's choice of delays, we can determine the start time of each activity in order to minimise this worst-case makespan.

Let us assume a fixed scenario $\theta\in\mathcal{U}(\Gamma)$ given by the vector $\delta\in[0,1]^{|V|}$. In this case, the inner minimisation problem becomes
\begin{align}
\min\ & S_{n+1}\\
\text{s.t. } & S_0 = 0\\
& S_j - S_i \geq \bar{\theta}_i + \delta_i\hat{\theta}_i - M(1-y_{ij}) & \forall (i,j)\in V^2\\
& S_i \ge 0 & \forall i\in V,
\end{align}
where $M$ is some number greater than or equal to the maximum possible minimum makespan. Taking the dual of this inner minimisation problem, we can find the following non-linear mixed-integer programming formulation for the adversarial sub-problem, first introduced in \cite{bruni2017adjustable}:
\begin{align}
\max\ & \sum_{(i,j)\in V^2}\left(\bar{\theta}_i+\delta_i\hat{\theta}_i - M(1-y_{ij})\right)\alpha_{ij}\\
\text{s.t. } & \sum_{(i,j)\in V^2}\alpha_{ij}-\sum_{(j,i)\in V^2}\alpha_{ji}=0 & \forall j \in V\\
& \sum_{(0,i)\in V^2}\alpha_{0i}=1\\
& \sum_{(i,n+1)\in V^2}\alpha_{i,n+1}=1\\
& \sum_{i\in V}\delta_i\leq \Gamma\\
& 0 \leq \delta_i \leq 1 & \forall i \in V\\
& \alpha_{ij}\in\{0,1\} &  \forall (i,j)\in V^2.
\end{align}
This can be viewed as a longest-path problem, where up to $\Gamma$ units of delay can be distributed among activities by the adversary.

As shown by \cite{bruni2017adjustable}, this model can be linearised as follows:
\begin{align}
\max\ & \sum_{(i,j)\in V^2}\left(\bar{\theta}_i\alpha_{ij}+\hat{\theta}_iw_{ij} - M(1-y_{ij})\alpha_{ij}\right)\label{eqn:linearised1}\\
\text{s.t. } & \sum_{(i,j)\in V^2}\alpha_{ij}-\sum_{(j,i)\in V^2}\alpha_{ji}=0 & \forall j \in V\label{eqn:linearised2}\\
& \sum_{(0,i)\in V^2}\alpha_{0i}=1\label{eqn:linearised3}\\
& \sum_{(i,n+1)\in V^2}\alpha_{i,n+1}=1\label{eqn:linearised4}\\
& w_{ij}\leq \delta_i & \forall (i,j) \in V^2\label{eqn:linearised5}\\
& w_{ij}\leq \alpha_{ij} & \forall (i,j) \in V^2\label{eqn:linearised6}\\
& \sum_{i\in V}\delta_i\leq \Gamma\label{eqn:linearised7}\\
& 0 \leq \delta_i \leq 1 & \forall i \in V\label{eqn:linearised8}\\
& \alpha_{ij}\in\{0,1\}& \forall (i,j)\in V^2\label{eqn:linearised9}\\
& w_{ij} \geq 0 & \forall (i,j) \in V^2\label{eqn:linearised10}.
\end{align}
It is claimed in Proposition 4 of \cite{bruni2017adjustable} that this problem is equivalent to its linear relaxation, where $\alpha_{ij}\in[0,1]$ for all $(i,j)\in V^2$. This, however, is not the case, as the following counter-example demonstrates.

\begin{figure}[htb]
	\centering
	\subfloat[$\alpha_{ij}\in\{0,1\}\quad\forall(i,j)\in V^2$ \label{fig:counter_example_a}]{
		\includegraphics[scale=0.7]{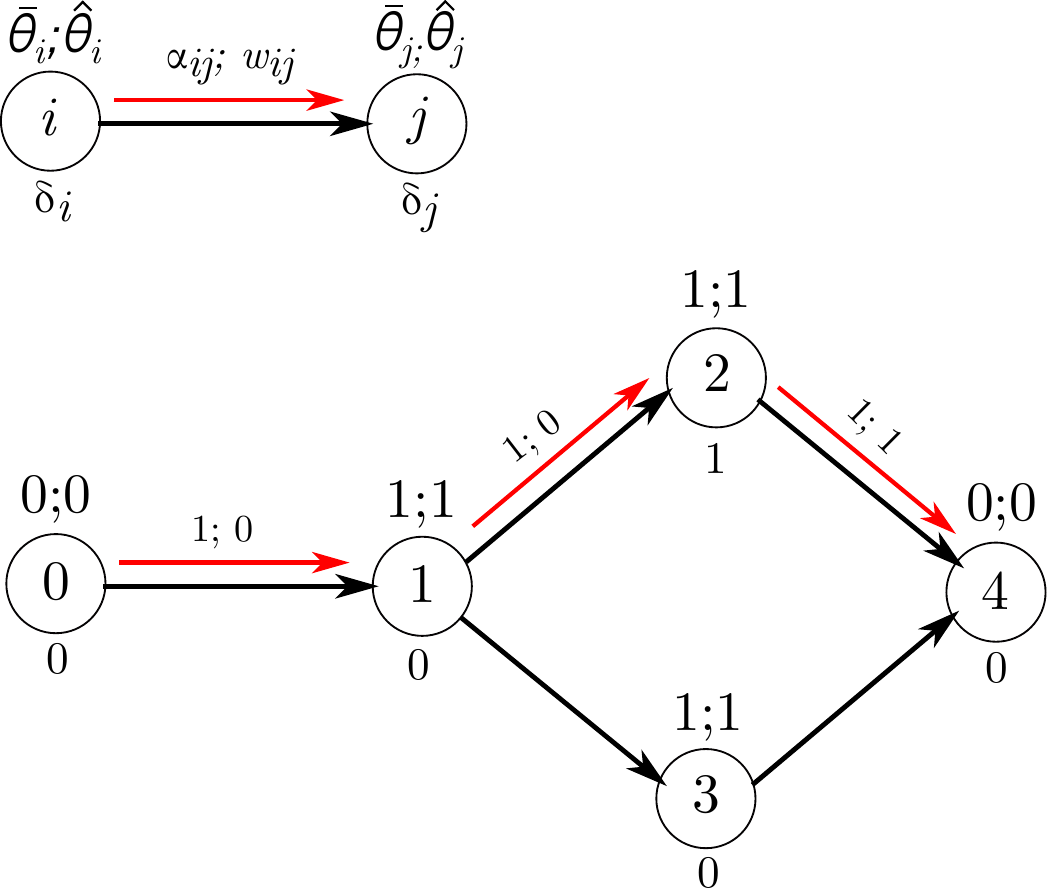}
	}

	\subfloat[{$\alpha_{ij}\in[0,1]\quad \forall (i,j)\in V^2$} \label{fig:counter_example_b}]{
		\includegraphics[scale=0.7]{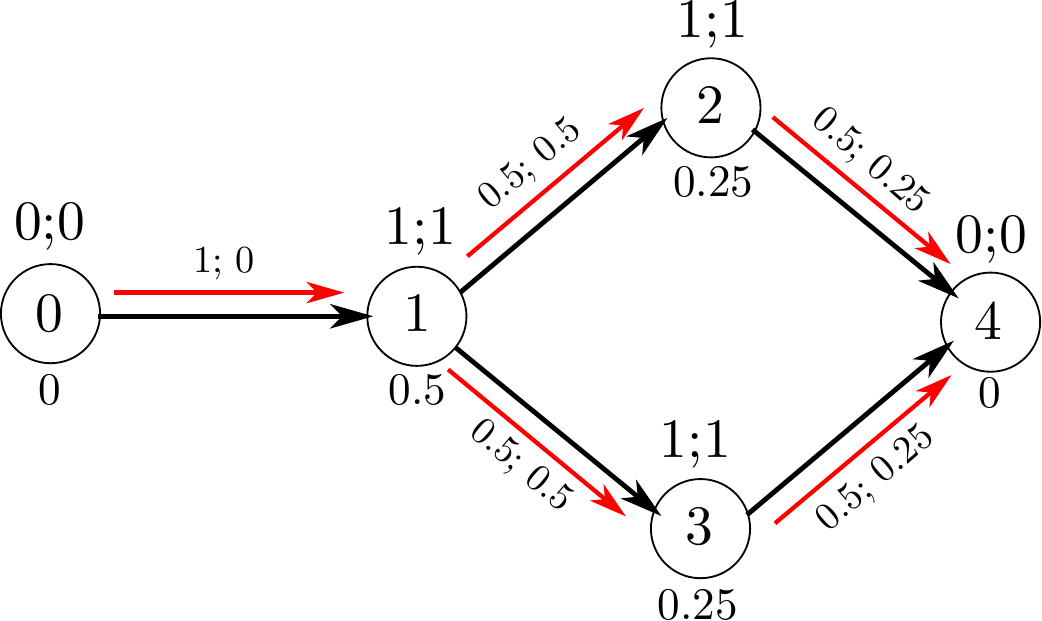}
	}
	\caption{Counter-example showing that model (\ref{eqn:linearised1})-(\ref{eqn:linearised10}) is not equivalent to its linear relaxation.}
	\label{fig:counter_example}
\end{figure}

Figure~\ref{fig:counter_example} shows a project with three non-dummy activities, each with a nominal duration of $\bar{\theta}_i=1$, and a maximum deviation of $\hat{\theta}_i=1,\,i=1,2,3$. Suppose a feasible first-stage solution has been found, resulting in the network shown in Figure \ref{fig:counter_example}. We consider this problem from the point of view of the adversary, who wishes to distribute up to $\Gamma=1$ units of delay, in order to maximise the minimum makespan. If (\ref{eqn:linearised1})-(\ref{eqn:linearised10}) is equivalent to its linear relaxation, then the adversary gains no advantage by choosing $\alpha\in(0,1)$ and splitting the unit flow on its route from the source-node 0 to the sink-node 4. However, as can be seen with this example, the adversary does in fact obtain an advantage.

In Figure~\ref{fig:counter_example_a}, $\alpha_{ij}\in \{0,1\}$ for each $(i,j)\in V^2$, and hence the adversary is limited to routing the unit flow through the network via a single path. A worst-case delay in this scenario is that the unit of available delay is entirely assigned to activity 2. Hence, $\delta_2=1$, whilst $\delta_1=\delta_3=0$. Minimising the worst-case makespan in this scenario, we get $(\bar{\theta}_1\alpha_{12}+\hat{\theta}_1w_{12})+(\bar{\theta}_2\alpha_{24}+\hat{\theta}_2w_{24})=(1+0)+(1+1)=3$.

In Figure~\ref{fig:counter_example_b}, $\alpha_{ij}\in[0,1]$ for each $(i,j)\in V^2$, and the adversary is able to split the unit flow into multiple fractional paths on its route through the network. In this case, the adversary can distribute the unit of delay so that $\delta_1=0.5$, $\delta_2=0.25$, and $\delta_3=0.25$. In this scenario, the minimum makespan is $(\bar{\theta}_1\alpha_{12}+\hat{\theta}_1w_{12}) + (\bar{\theta}_1\alpha_{13}+\hat{\theta}_1w_{13})+(\bar{\theta}_2\alpha_{24}+\hat{\theta}_2w_{24})+(\bar{\theta}_3\alpha_{34}+\hat{\theta}_3w_{34})=(0.5+0.5)+(0.5+0.5)+(0.5+0.25)+(0.5+0.25)=3.5$, showing that problem (\ref{eqn:linearised1})-(\ref{eqn:linearised10}) is not equivalent to its linear relaxation.

Note that \cite{bruni2017adjustable} attempt to prove that model (\ref{eqn:linearised1})-(\ref{eqn:linearised10}) is equivalent to its linear relaxation, and therefore polynomially solvable, by showing that the corresponding constraint matrix is totally unimodular. In Appendix~\ref{sec:appendix}, we identify an error with this proof and show that the constraint matrix is not totally unimodular. This result is consistent with the above counter-example.

Since problem (\ref{eqn:linearised1})-(\ref{eqn:linearised10}) is not equivalent to its linear relaxation, we cannot apply strong-duality to get an equivalent minimisation problem. Therefore, in order to obtain a compact reformulation of the two-stage robust RCPSP, an alternative reformulation of the adversarial sub-problem is required.

A dynamic programming procedure for solving problem (\ref{eqn:linearised1})-(\ref{eqn:linearised10}) when $\Gamma\in \mathbb{Z}$ is presented in \cite{bruni2017adjustable}. This procedure works by considering $\Gamma+1$ paths from the source node 0 to node $i$, for each $i\in V$, where each path $\pi_i^\gamma,\,\gamma=0,\dots,\Gamma$, is characterised by the inclusion of exactly $\gamma$ delayed activities. Given a path $\pi_i^\gamma$, its extension to each successor node $j\in Succ_i$ is evaluated by considering two possibilities: either the successor activity $j$ is delayed, resulting in the path $\pi_j^{\gamma+1}$, or it is not delayed, resulting in the path $\pi_j^\gamma$. Hence, the dynamic programming algorithm has a state $ST(j,\gamma)$ for each node $j$ at level $\gamma$, and the value of each state $V(ST(j,\gamma))$ is computed through the following recursion:
\begin{align}
	&V(ST(0,0))=0,\\
	\begin{split}
	V(ST(j,\gamma))=\max_{i:(i,j)\in E\cup X}\Big\{\max\Big(V(ST(i,\gamma)),\,V(ST(&i,\gamma-1))+\bar{\theta}_i+\hat{\theta}_i\Big)\Big\},\\
	&\quad \forall j\in V\setminus\{0\},\,\gamma=1,\dots,\Gamma\\
	\end{split}\\
	&V(ST(j,0))=\max_{i:(i,j)\in E\cup X}\Big\{V(ST(i,0)) + \bar{\theta}_i\Big\}.
\end{align}

\begin{figure}[htbp]
	\centering
	\includegraphics[scale=0.7]{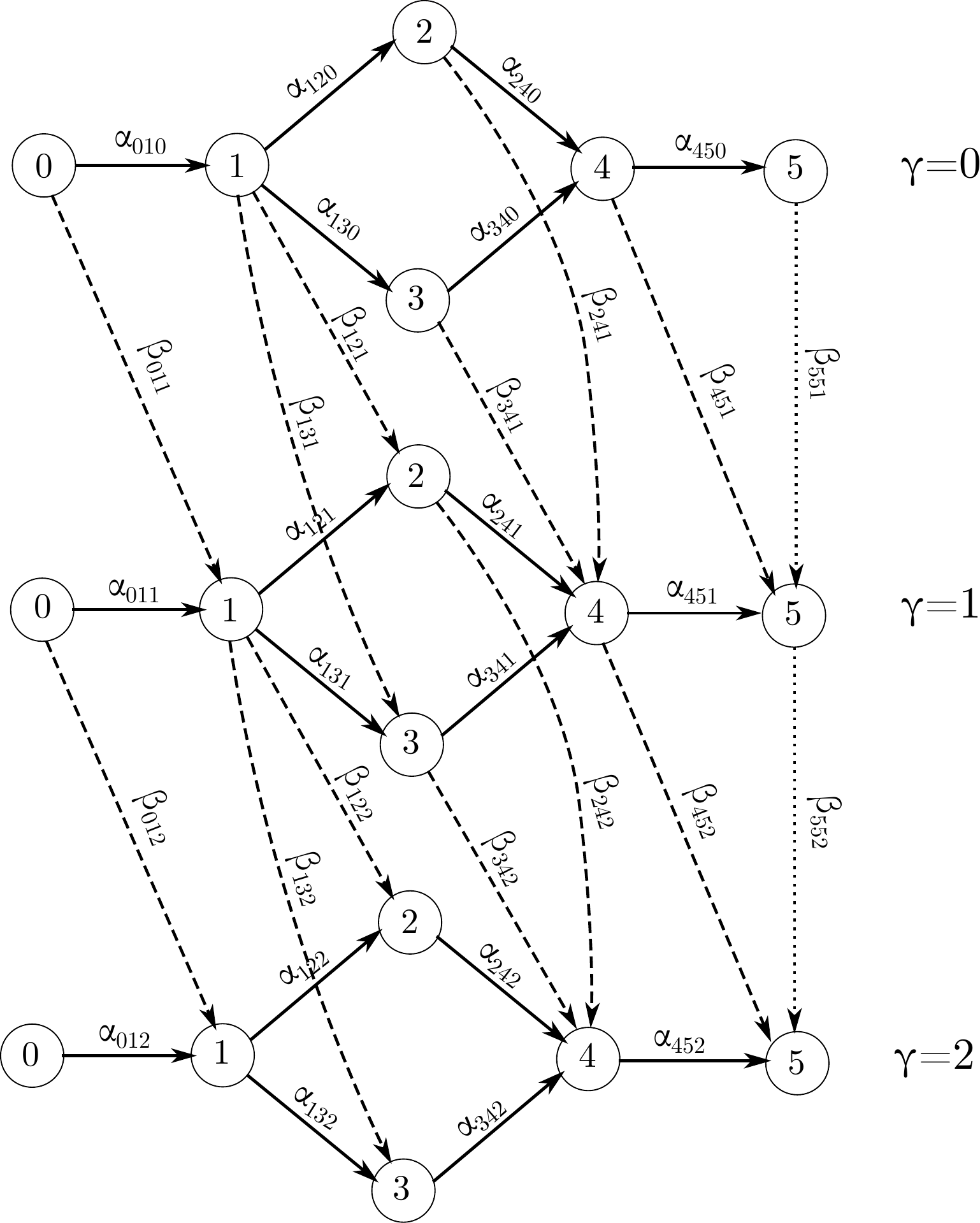}
	\caption{Example augmented graph for a project with four non-dummy activities, and where up to $\Gamma=2$ activities can reach their worst-case durations.}
	\label{fig:aux_graph}
\end{figure}

This dynamic programming algorithm can be viewed as finding the critical path through the augmented project network built from $\Gamma+1$ copies of the original project network (an example of such a network is shown in Figure~\ref{fig:aux_graph}). The inclusion of an inter-level arc, e.g. a dashed arc in Figure~\ref{fig:aux_graph}, in the critical path corresponds to the delay of the activity at the origin of that arc.

Since the second stage problem is simply a longest-path problem on this augmented network, it can be recast into the following mixed-integer linear program:
\begin{align}
\max & \sum_{(i,j)\in V^2}\sum_{\gamma=0}^{\Gamma}(\bar{\theta}_i-M(1-y_{ij}))\alpha_{ij\gamma}+\sum_{(i,j)\in V^2}\sum_{\gamma=1}^{\Gamma}(\bar{\theta}_i+\hat{\theta}_i-M(1-y_{ij}))\beta_{ij\gamma} \hspace*{-5cm} \label{eqn:sub_reformulation1}\\
\text{s.t. } &
% \begin{split}
    \sum_{(j,i)\in V^2}\alpha_{ji\gamma} + \sum_{(j,i)\in V^2}\beta_{ji,\gamma+1} - \sum_{(i,j)\in V^2}\alpha_{ij\gamma} - \sum_{(i,j)\in V^2}\beta_{ij\gamma} = 0 \hspace*{-2cm}\nonumber \\
& & \hspace*{-7cm}\forall j\in V,\,\gamma=1,\dots,\Gamma-1 \label{conflow}\\
% &\nu_{\gamma+1} - \sum_{(i,n+1)\in V^2}\alpha_{i,n+1,\gamma} - \sum_{(i,n+1)\in V^2}\beta_{i,n+1,\gamma}=0 & \gamma =1,\dots,\Gamma-1 \label{conflow_n+1}\\
% \end{split}\\
& \sum_{(j,i)\in V^2}\alpha_{ji0} + \sum_{(j,i)\in V^2}\beta_{ji1} - \sum_{(i,j)\in V^2}\alpha_{ij0} = 0 & \forall j\in V \label{conflow2}\\
& \sum_{(j,i)\in V^2}\alpha_{ji\Gamma} - \sum_{(i,j)\in V^2}\alpha_{ij\Gamma} - \sum_{(i,j)\in V^2}\beta_{ij\Gamma} = 0 & \forall j\in V \label{conflow3}\\
& \sum_{(0,i)\in V^2}\alpha_{0i0} + \sum_{(0,i)\in V^2}\beta_{0i1} = 1 \label{constart}\\
& \sum_{(i,n+1)\in V^2}\alpha_{i,n+1,\Gamma} + \sum_{(i,n+1)\in V^2}\beta_{i,n+1,\Gamma} = 1 \label{conend} \hspace*{-7cm}\\
& \alpha_{ij\gamma}\in\{0,1\} & \forall (i,j)\in V^2,\, \gamma =0,\dots,\Gamma\\
& \beta_{ij\gamma}\in\{0,1\} & \forall (i,j)\in V^2,\,\gamma=1,\dots,\Gamma\label{eqn:sub_reformulation_end}
\end{align} 
where $\alpha_{ij\gamma}$ is the flow from node $i$ to node $j$ in level $\gamma$ and $\beta_{ij\gamma}$ is the flow from node $i$ in level $\gamma-1$ to node $j$ in level $\gamma$.
The constraints model a unit flow through the augmented network from node 0 in level 0 (Constraint (\ref{constart})) to node $n+1$ in level $\Gamma$ (Constraint (\ref{conend})). Constraints (\ref{conflow})
% and (\ref{conflow_n+1}) 
are flow-conservation constraints that ensure that for node each in level $\gamma=1,\dots,\Gamma-1$, the incoming flow from levels $\gamma$ and $\gamma-1$ must be equal to the outgoing flow to levels $\gamma$ and $\gamma+1$. Constraints (\ref{conflow2}) and (\ref{conflow3}) conserve flow over the nodes in the special cases of the first and last level, respectively.

Note that this model includes more $\alpha_{ij\gamma}$ and $\beta_{ij\gamma}$ variables than indicated in Figure~\ref{fig:aux_graph}, with the edges shown in Figure \ref{fig:aux_graph} corresponding to the edges for which $y_{ij}=1$. The edges that are not shown are penalised by constant $M$ in the objective~(\ref{eqn:sub_reformulation1}) when $y_{ij} = 0$. To ensure that it is always possible to find a path from node 0 in level 0 to node $n+1$ in level $\Gamma$ in the augmented network (if $\Gamma$ is larger than the number of activities included in the longest path from node 0 to node $n+1$ in the original project network, such a path may not be possible), the final sink nodes of each layer are connected by enforcing $y_{n+1,n+1}=1$ (see dotted arcs in Figure \ref{fig:aux_graph}). Since $\bar{\theta}_{n+1}+\hat{\theta}_{n+1} = 0$ these additional edges can be traversed at no extra cost to reach node $n+1$ in level $\Gamma$.

% In the case where $\Gamma$ is larger than the number of activities included in the longest path from node 0 to node $n+1$ in the original project network, it may be impossible to find a path through through the augmented network from node 0 in level 0 to node $n+1$ in level $\Gamma$. To avoid this pathological situation, arcs of length 0 have been inserted between the dummy sink nodes of levels 1 to $\Gamma$ of the augmented network (see dotted arc in Figure \ref{fig:aux_graph}). Now, if the longest path through the augmented network finishes at the sink node of a level before $\Gamma$ (i.e. not all the available delay was required in the worst-case scenario), the path can traverse these arcs at no extra cost to reach node $n+1$ in level $\Gamma$. Note that an arc is never required between the sink nodes of levels 0 and 1.

\subsection{Compact reformulation}\label{subsection:compact_reformulation}

Since problem (\ref{eqn:sub_reformulation1})-(\ref{eqn:sub_reformulation_end}) is simply a longest-path problem over an augmented project graph, it is equivalent to its linear relaxation where $\alpha_{ij\gamma}\in [0,1]$ for all $(i,j)\in V^2,\,\gamma = 0,\dots,\Gamma$, and $\beta_{ij\gamma}\in [0,1]$ for all $(i,j)\in V^2,\,\gamma = 1,\dots,\Gamma$. Hence, we can take the dual of this problem to get an equivalent formulation for the second-stage problem as a minimisation problem. 

The first-stage problem aims determine a sufficient selection $X\in \mathcal{X}$ that minimises the objective value of the second-stage objective value. This first-stage problem can be modelled with a flow-based formulation, as proposed by \cite{artigues2003insertion}. This formulation makes use of continuous resource flow variables $f_{ijk}$, which determine the amount of resource type $k\in K$ that is transferred upon the completion of activity $i$ to activity $j$. Additionally, binary variables $y_{ij}$ capture the choice of sufficient selection by representing precedence relationships of the extended project network.

Thus, having dualised the second-stage problem (\ref{eqn:sub_reformulation1})-(\ref{eqn:sub_reformulation_end}) into a minimisation problem, the first and second-stages can be combined to obtain the following compact reformulation of the full two-stage robust RCPSP with budgeted uncertainty:
\begin{align}
\min \ & S_{n+1,\Gamma}\label{eqn:compact1}\\
\text{s.t. } & S_{00}=0\label{eqn:compact2}\\
& S_{j\gamma} - S_{i\gamma} \geq \bar{\theta}_i - M(1-y_{ij}) & \forall (i,j)\in V^2,\,\gamma=0,\dots,\Gamma\label{eqn:compact_bigM1}\\
& S_{j,\gamma+1} - S_{i\gamma} \geq \bar{\theta}_i + \hat{\theta}_i - M(1-y_{ij}) & \forall (i,j)\in V^2,\,\gamma = 0,\dots,\Gamma-1\label{eqn:compact_bigM2}\\
% & S_{n+1,\gamma+1}-S_{n+1,\gamma}\geq 0 & \gamma = 1,\dots,\Gamma-1\label{eqn:compact_add}\\
& y_{ij}=1 & \forall (i,j) \in E\cup\{(n+1,n+1)\}\label{eqn:compact5}\\
& f_{ijk}\leq N_ky_{ij} & \forall (i,j)\in V^2,\, \forall k \in K\label{eqn:compact6}\\
&\sum_{i\in V}f_{ijk}=r_{jk} & \forall j \in V ,\, \forall k\in K\label{eqn:compact7}\\
&\sum_{j\in V}f_{ijk}=r_{ik} & \forall i \in V,\,\forall k\in K\label{eqn:compact8}\\
&S_{i\gamma}\geq 0 &\forall i\in V,\,\gamma \in 0,\dots, \Gamma\label{eqn:compact9}\\
&f_{ijk}\geq 0 & \forall (i,j)\in V^2 ,\, \forall k \in K\label{eqn:compact10}\\
&y_{ij}\in\{0,1\} & \forall (i,j)\in V^2\label{eqn:compact11},
\end{align}
where $M$, as before, is chosen to be greater than or equal to the maximum possible minimum makespan, and $N_k$ is some number greater than or equal to $R_k$. Constraints~(\ref{eqn:compact2})-(\ref{eqn:compact_bigM2}) are the dual constraints of the second-stage problem (\ref{eqn:sub_reformulation1})-(\ref{eqn:sub_reformulation_end}), and ensure that activity start time respect the project precedence relationships. Constraints~(\ref{eqn:compact5}) capture the original project precedences, whilst constraints~(\ref{eqn:compact6})-(\ref{eqn:compact8}) are resource flow constraints. Constraints~(\ref{eqn:compact6}) ensure that resource flow respects precedence relationships, and constraints~(\ref{eqn:compact7}) and (\ref{eqn:compact8}) conserve flow into and out of each node, respectively.

It is important to note that this basic model does not enforce the transitivity of the $y$-variables. Instead, the model captures the extended project network in terms of the $y$-variables with constraints (\ref{eqn:compact5}) and (\ref{eqn:compact6}), and ensures the feasibility of activity start-times with respect to this extended network through constraints (\ref{eqn:compact_bigM1}) and (\ref{eqn:compact_bigM2}). In Section \ref{section:computational_experiments} the computational benefits of extending model (\ref{eqn:compact1})-(\ref{eqn:compact11}) to include explicit transitivity constraints on the $y$-variables is examined.

\section{Computational experiments} \label{section:computational_experiments}

This section compares results obtained by solving the compact robust counterpart (\ref{eqn:compact1})-(\ref{eqn:compact11}), and three slight extensions to this model, with the current state-of-the-art approach to solving the two-stage robust RCPSP proposed in \cite{bruni2018computational}. Before outlining the proposed extensions to the basic model detailed in the previous section, we introduce the test instances used in this computational study. 

\subsection{Instances}

The test instances used in this computational study have been converted from deterministic RCPSP instances involving 30 activities, taken from the PSPLIB (\cite{kolisch1997psplib}, \url{http://www.om-db.wi.tum.de/psplib/}). The difficulty of these instances is measured and controlled by the following three parameters:

\begin{enumerate}
	\item Network complexity $NC\in\{1.5, 1.8, 2.1\}$. This measures the average number of non-redundant (i.e. non-transitive) arcs per activity.
	\item Resource factor $RF\in\{0.25,0.5,0.75,1\}$. This measures the average proportion of resource types for which a non-dummy activity has a non-zero requirement.
	\item Resource strength $RS\in\{0.2,0.5,0.7,1\}$. This measures the restrictiveness of the availability of the resources, with a smaller $RS$ value indicating a more constrained project instance.
\end{enumerate}
The PSPLIB contains a set of 10 instances for each of the 48 possible combinations of instance parameters. 

The maximum deviation of the duration of each activity is set to be $\hat{\theta}=\bigl\lceil\bar{\theta}/2\bigr\rceil$. For each of the 480 deterministic RCPSP instances in the PSPLIB, three robust counterparts have been generated by considering $\Gamma \in \{3,5,7\}$, resulting in a total of 1440 test instances. The sets of 30 robust counterparts for each combination of instance parameters are labelled J301, J302, $\dots$, J3048. Note that the instances used in this computational study are identical to the instances used in \cite{bruni2017adjustable} and \cite{bruni2018computational}.    

\subsection{Implementations}\label{section:implementations}

The following section compares the performance of model (\ref{eqn:compact1})-(\ref{eqn:compact11}) with that of three slight extensions to this model. Here, we outline these extensions and clarify details regarding the practical implementation of these models.

The first variant of the basic model (\ref{eqn:compact1})-(\ref{eqn:compact11}) includes the following transitivity constraints on the $y$-variables:
\begin{align}
	& y_{ij}+y_{ji}\leq 1 &\forall (i,j)\in V^2\setminus\{(n+1,n+1)\}\label{eqn:trans1}\\
	& y_{ij}\geq y_{il} + y_{lj}-1  &\forall (i,l,j)\in V^3\label{eqn:trans2}.
\end{align}
As explained in Section \ref{subsection:compact_reformulation}, these transitivity of the $y$-variables is not strictly necessary to ensure the feasibility of the activity start-times. We include them as an extension to model (\ref{eqn:compact1})-(\ref{eqn:compact11}) in order to assess their impact on the computational performance of the model.

The second extension involves the provision of a heuristic warm-start solution to the solver software. This heuristic solution is obtained with the following procedure:
\begin{enumerate}
	\item Given an uncertain RCPSP instance, a heuristic solution is found to the corresponding deterministic instance using the latest-finish-time (LFT) priority-rule heuristic \citep{kolisch1996serial}. 
	\item From this solution, a feasible set of $y$-variables is obtained by setting$$y_{ij}=
\begin{cases}
	1 \quad\textnormal{if } s_j \geq f_i\\
	0 \quad\textnormal{otherwise},
\end{cases}
$$ 
where $s_j$ is the start time of activity $j$, and $f_i$ is the finish time of activity $i$.
	\item These $y$ variables are passed to the basic model (\ref{eqn:compact1})-(\ref{eqn:compact11}), which is solved to provide a feasible warm-start solution.
\end{enumerate}

This warm-start solution can be used to tighten the big-$M$ constraints (\ref{eqn:compact_bigM1}) and (\ref{eqn:compact_bigM2}), and thereby further improve the basic model. This is achieved by setting $M_{ij}=LF^*_i-ES_j$ for each $(i,j)\in V^2$, where $ES_j$ is the earliest start time of activity $j$, and $LF^*_i$ is the latest finish time of activity $i$, calculated relative to the makespan of the warm-start solution. These values are computed recursively via a forward-pass and backward-pass of the project network, respectively.

Note that, although the $S$-variables of model (\ref{eqn:compact1})-(\ref{eqn:compact11}) are in general continuous, for the purposes of this computational study, the $S$-variables have been set to be integer. Since $\hat{\theta}=\bigl\lceil\bar{\theta}/2\bigr\rceil\in \mathbb{Z}$ for the instances solved in this study, the correctness of the model is unaffected by this specification.

In summary, the following section presents results from the following five solution approaches:
\begin{enumerate}
	\item Basic model (\ref{eqn:compact1})-(\ref{eqn:compact11}),
	\item Basic model with transitivity constraints, i.e (\ref{eqn:compact1})-(\ref{eqn:trans2}),
	\item Basic model with warm-start,
	\item Basic model with warm-start and transitivity constraints,
	\item Primal method from \cite{bruni2018computational}. This is the strongest existing approach for solving the two-stage robust RCPSP.
\end{enumerate}

All the models proposed in this paper have been solved using Gurobi 9.0.1, running on 4 cores of a 2.30GHz Intel Xeon CPU, limited to 16GB RAM. Note that the specifications of this machine have been chosen to be as similar as possible to that of the CPU used in the experiments performed in \cite{bruni2017adjustable} and \cite{bruni2018computational}. A limit of 20 minutes was imposed on the solution time of each model, the same as used for the experiments performed in \cite{bruni2017adjustable} and \cite{bruni2018computational}. Results for the primal method have been reproduced from \cite{bruni2018computational}. 

\subsection{Results}

In this section, we first present and analyse results from solving model (\ref{eqn:compact1})-(\ref{eqn:compact11}) and the three variants proposed in the previous section, before we compare these results with those from the current best iterative algorithm presented in \cite{bruni2018computational}. 

\begin{figure}[h]
	\centering
	\includegraphics[scale=0.75]{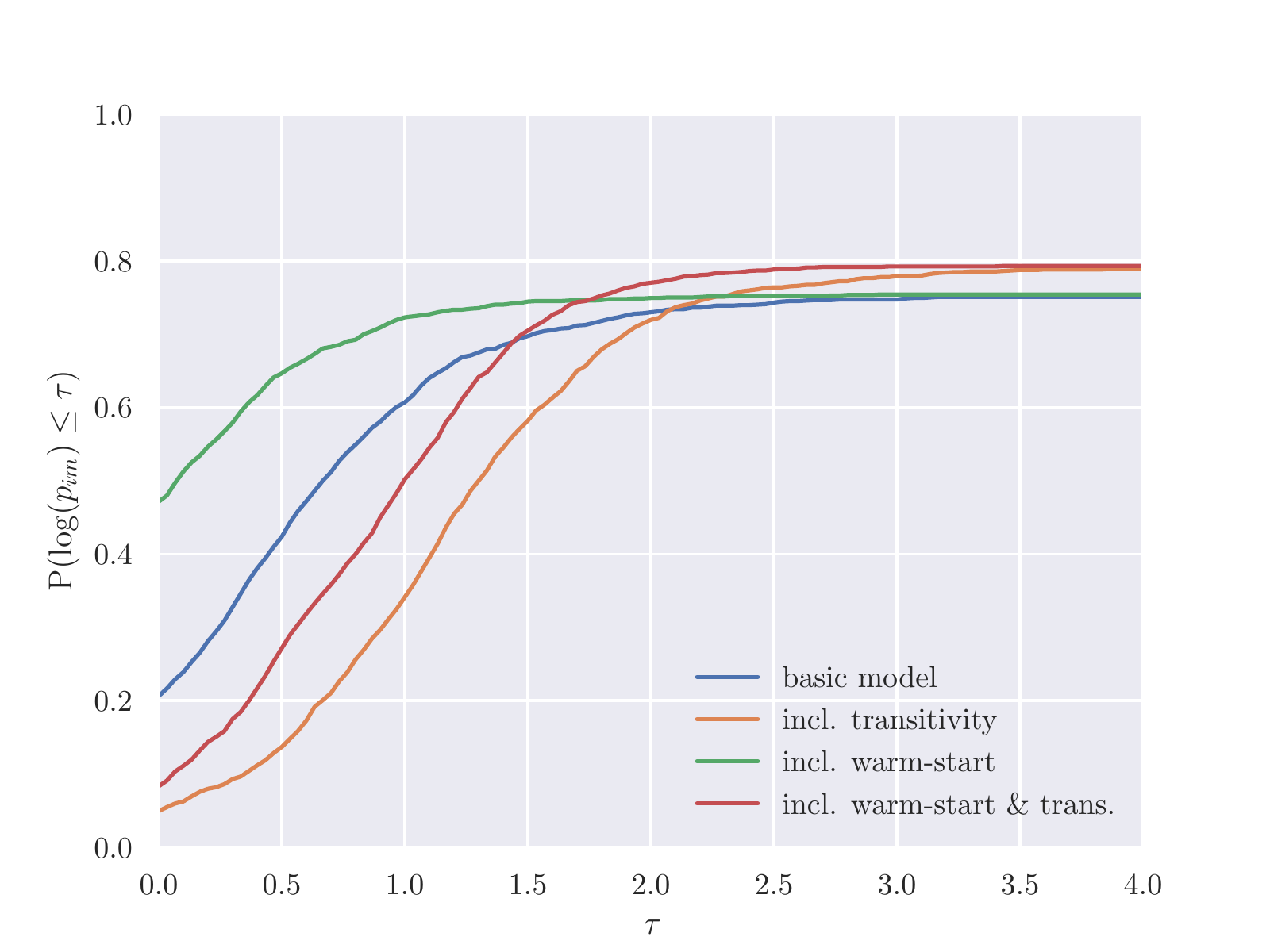}
	\caption{Performance profile of relative solution times.}
	\label{fig:pp_time_ratios}
\end{figure}

We start by considering the performance profile \citep{dolan2002benchmarking} plot shown in Figure \ref{fig:pp_time_ratios}. The performance profile uses the \textit{performance ratio} as a measure by which the different models can be compared. The performance ratio of model $m\in\mathcal{M}$ for problem instance $i\in\mathcal{I}$ is defined to be
$$p_{im}=\frac{t_{im}}{\min_{m\in \mathcal{M}}t_{im}},$$
where $t_{im}$ is the time required to solve instance $i$ using model $m$. If model $m$ is unable to solve instance $i$ to optimality within the 20 minute time-limit, then $p_{im}=P$, where $P\geq \max_{i,m} r_{im}$. The performance profile of model $m\in \mathcal{M}$ is defined to be the function $$\rho_m(\tau)=\frac{|\{p_{im}\leq \tau\,:\, i\in \mathcal{I}\}|}{|\mathcal{I}|},$$ i.e. the probability that the performance ratio of model $m$ is within a factor $\tau$ of the best performance ratio. The performance profile in Figure \ref{fig:pp_time_ratios} has been plotted on the log scale for clarity. 
 
It is clear from Figure~\ref{fig:pp_time_ratios} that the provision of a heuristic warm-start solution improves solution time, with the models that make use of a warm-start solution being faster to solve for a greater proportion of instances that their respective models without a warm-start. It can also be seen that the models that make use of transitivity constraints are slower to solve to optimality for a greater proportion of instances than their respective models that do not use transitivity constraints. However, the inclusion of transitivity constraints does increase the proportion of instances that can be solved to optimality, by 5.3\% for the basic model, and by 5.2\% for the model with warm-start.  

\begin{figure}[h]
	\centering
	\includegraphics[scale=0.75]{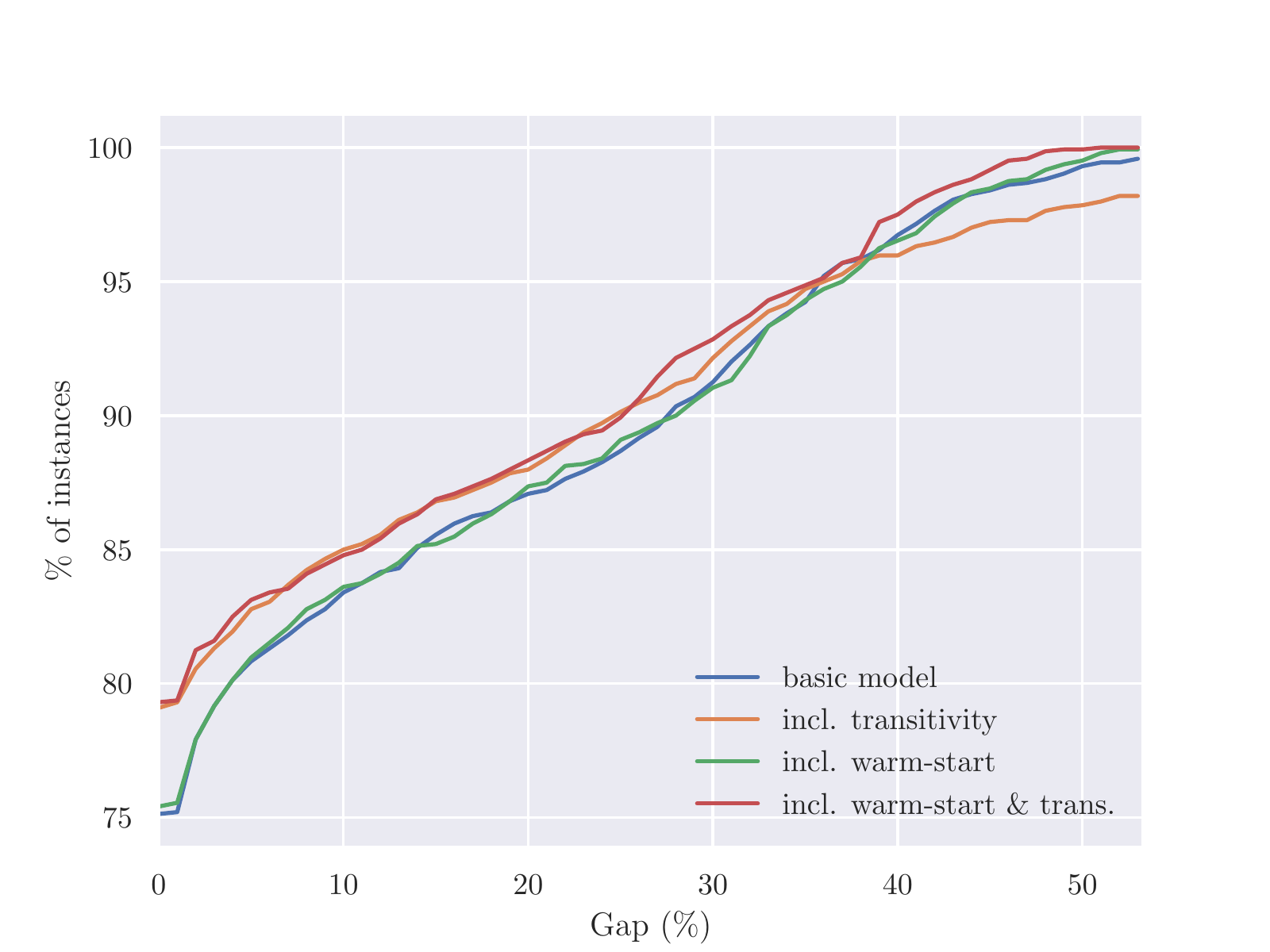}
	\caption{Cumulative percentage of instances solved to within given gap of optimality within time-limit.}
	\label{fig:pp_gap}
\end{figure}

Figure~\ref{fig:pp_gap} plots the cumulative percentage of instances solved to within a given optimality gap within the 20 minute time-limit. Note that the left-hand y-intercept of this figure gives the same information as the right-hand y-intercept in Figure~\ref{fig:pp_time_ratios}, that is, the proportion of instances solved to optimality using each model. Looking at Figure \ref{fig:pp_gap}, it can be seen that as well as increasing the proportion of instances that can be solved to optimality, the inclusion of transitivity constraints increases the proportion of instances that can be solved to within a given optimality gap. Of the 255 instances for which an optimal solution was unable to be found with any model, but for which a feasible solution was found using all models, the average optimality gap was 24.53\% for the basic model, 22.80\% with the inclusion of transitivity constraints, 24.71\% with the inclusion of a warm-start solution, and 22.36\% with the inclusion of both a warm-start solution and transitivity constraints. Note however that the basic model fails to find a feasible solution for 3 instances, whilst the model that includes transitivity constraints only fails to find a feasible solution for 24 instances. The other two variants find feasible solutions to all 1440 instances.

From Figures \ref{fig:pp_time_ratios} and \ref{fig:pp_gap}, we can see that the inclusion of a warm-start solution and transitivity constraints in model (\ref{eqn:compact1})-(\ref{eqn:compact11}), is the best performing variant: it solves the greatest number of instances to optimality, is the strongest performing model over the instances which no model can solve to optimality, and is significantly faster to solve than the transitive model without a warm-start. 

Of the 1440 test instances, 1160 have been solved to optimality within the time-limit by at least one of the four variants of model (\ref{eqn:compact1})-(\ref{eqn:compact11}) proposed in this paper. The strongest single model is the model which includes transitivity constraints and makes use of a heuristic warm-start solution, solving 1142 instances to optimality.

In Table \ref{table:results_comparison2}, we now compare the performance of the basic model (\ref{eqn:compact1})-(\ref{eqn:compact11}) and its strongest extension, with the results of the strongest existing algorithm for the two-stage robust RCPSP, the \textit{primal method} \citep{bruni2018computational}. For each set of test instances, J301, $\dots$, J3048, Table \ref{table:results_comparison2} reports instance parameters (NC, RF, RS), the average CPU time required to solve the instances that were solved to optimality ($time$), the average optimality gap for those instances which were not solved to optimality but for which a feasible solution was obtained ($gap$), and the number of instances solved to optimality ($\#solv$).

Note that in four of the most challenging instance sets, J3013, J3021, J3029, J3041, the results in Table \ref{table:results_comparison2} show that the primal method solves one or two instances to optimality, sometimes outperforming the model proposed in this paper over these instance sets. However, these optimal solutions are obtained whilst simultaneously reaching the maximum time-limit of 1200 seconds, and it is therefore unclear whether or not this is a numerical inaccuracy in the results presented in \cite{bruni2018computational}. 

The results in this table show that the models proposed in this paper solve almost 50\% more instances than the primal method, and do so in a considerably shorter computation time. These results confirm the strength of the new model proposed in this paper.

\newpage
\newgeometry{left=0.5in, right=0.5in, top=1in}

\begin{table}[]
\centering
\begin{tabular}{llllrrrrrrrrrrr}
	\hline\hline
      &     &      &     & \multicolumn{3}{p{3.5cm}}{\centering Primal Method \\ \citep{bruni2018computational}} && \multicolumn{3}{p{3.5cm}}{\centering\ \\[-1.5mm]  Basic model (\ref{eqn:compact1})-(\ref{eqn:compact11})} && \multicolumn{3}{p{3.6cm}}{\centering \ \\[-1.5mm] {\small incl. warm-start + trans.}} \\
\cline{5-7} \cline{9-11} \cline{13-15}
      & $NC$  & $RF$   & $RS$  & $time$      & $gap$      & $\#solv$   && $time$     & $gap$      & $\#solv$  && $time$         & $gap$          & $\#solv$       \\
      \hline\hline
j301  & 1.50 & 0.25 & 0.20 & 196.9     & 5.54     & 21         && 6.96     &          & 30        && 19.69        &              & 30             \\
j302  & 1.50 & 0.25 & 0.50 & 120.42    & 3.64     & 28         && 2.68     &          & 30        && 6.57         &              & 30             \\
j303  & 1.50 & 0.25 & 0.70 & 13.04     & 4.55     & 29         && 1.11     &          & 30        && 2.42         &              & 30             \\
j304  & 1.50 & 0.25 & 1.00 & 5.52      & 11.80    & 27         && 0.78     &          & 30        && 1.62         &              & 30             \\
j305  & 1.50 & 0.50 & 0.20 &           & 15.92    & 0          && 676.70   & 16.88    & 1         && 822.10       & 18.12        & 8              \\
j306  & 1.50 & 0.50 & 0.50 & 358.61    & 12.84    & 3          && 119.58   & 1.52     & 29        && 102.27       &              & 30             \\
j307  & 1.50 & 0.50 & 0.70 & 260.77    & 8.20     & 19         && 8.10     &          & 30        && 10.74        &              & 30             \\
j308  & 1.50 & 0.50 & 1.00 & 59.06     & 6.93     & 22         && 1.45     &          & 30        && 1.95         &              & 30             \\
j309  & 1.50 & 0.75 & 0.20 &           & 10.19    & 0          &&          & 33.59    & 0         &&              & 30.71        & 0              \\
j3010 & 1.50 & 0.75 & 0.50 & 313.09    & 22.19    & 2          && 257.55   & 3.59     & 13        && 369.89       & 5.58         & 20             \\
j3011 & 1.50 & 0.75 & 0.70 & 232.82    & 12.83    & 7          && 94.28    & 1.46     & 28        && 54.21        & 1.38         & 28             \\
j3012 & 1.50 & 0.75 & 1.00 & 129.84    & 4.90     & 26         && 1.89     &          & 30        && 2.38         &              & 30             \\
j3013 & 1.50 & 1.00 & 0.20 & 1200      & 52.29    & 1          &&          & 38.02    & 0         &&              & 37.27        & 0              \\
j3014 & 1.50 & 1.00 & 0.50 & 440.91    & 20.29    & 2          && 293.75   & 5.74     & 7         && 333.90       & 5.73         & 12             \\
j3015 & 1.50 & 1.00 & 0.70 & 334.16    & 9.33     & 12         && 10.16    & 3.69     & 27        && 13.41        & 3.72         & 27             \\
j3016 & 1.50 & 1.00 & 1.00 & 97.46     & 7.43     & 27         && 1.33     &          & 30        && 2.79         &              & 30             \\
j3017 & 1.80 & 0.25 & 0.20 & 157.88    & 2.29     & 28         && 4.20     &          & 30        && 7.47         &              & 30             \\
j3018 & 1.80 & 0.25 & 0.50 & 18.26     &          & 30         && 1.31     &          & 30        && 2.20         &              & 30             \\
j3019 & 1.80 & 0.25 & 0.70 & 26.67     & 10.53    & 29         && 0.80     &          & 30        && 1.59         &              & 30             \\
j3020 & 1.80 & 0.25 & 1.00 & 8.23      & 5.75     & 28         && 0.40     &          & 30        && 1.28         &              & 30             \\
j3021 & 1.80 & 0.50 & 0.20 & 1200      & 9.94     & 2          && 503.17   & 11.65    & 10        && 462.90       & 12.59        & 18             \\
j3022 & 1.80 & 0.50 & 0.50 & 232.52    & 10.66    & 10         && 45.13    &          & 30        && 43.52        &              & 30             \\
j3023 & 1.80 & 0.50 & 0.70 & 145.11    & 4.29     & 24         && 2.71     &          & 30        && 4.65         &              & 30             \\
j3024 & 1.80 & 0.50 & 1.00 & 48.19     & 8.51     & 26         && 0.95     &          & 30        && 1.75         &              & 30             \\
j3025 & 1.80 & 0.75 & 0.20 &           & 13.15    & 0          &&          & 31.71    & 0         &&              & 29.77        & 0              \\
j3026 & 1.80 & 0.75 & 0.50 & 490.81    & 9.49     & 9          && 128.19   & 2.91     & 26        && 119.50       & 1.54         & 29             \\
j3027 & 1.80 & 0.75 & 0.70 & 128.64    & 7.52     & 16         && 3.29     &          & 30        && 4.15         &              & 30             \\
j3028 & 1.80 & 0.75 & 1.00 & 63.92     & 6.14     & 27         && 0.91     &          & 30        && 1.28         &              & 30             \\
j3029 & 1.80 & 1.00 & 0.20 & 1200      & 10.86    & 1          &&          & 40.12    & 0         &&              & 39.23        & 0              \\
j3030 & 1.80 & 1.00 & 0.50 &           & 19.98    & 0          && 785.57   & 5.01     & 3         && 774.16       & 5.08         & 8              \\
j3031 & 1.80 & 1.00 & 0.70 & 87.19     & 11.34    & 9          && 6.41     & 4.79     & 24        && 43.75        & 3.81         & 25             \\
j3032 & 1.80 & 1.00 & 1.00 & 45.52     & 10.99    & 26         && 1.00     &          & 30        && 1.16         &              & 30             \\
j3033 & 2.10 & 0.25 & 0.20 & 28.35     &          & 30         && 1.58     &          & 30        && 2.01         &              & 30             \\
j3034 & 2.10 & 0.25 & 0.50 & 10.37     & 2.35     & 29         && 0.66     &          & 30        && 0.79         &              & 30             \\
j3035 & 2.10 & 0.25 & 0.70 & 27.65     & 10.98    & 27         && 0.54     &          & 30        && 0.65         &              & 30             \\
j3036 & 2.10 & 0.25 & 1.00 & 20.52     &          & 30         && 0.29     &          & 30        && 0.44         &              & 30             \\
j3037 & 2.10 & 0.50 & 0.20 & 906.85    & 7.29     & 7          && 256.89   & 16.73    & 18        && 317.67       & 9.89         & 23             \\
j3038 & 2.10 & 0.50 & 0.50 & 239.89    & 6.81     & 23         && 11.63    &          & 30        && 12.53        &              & 30             \\
j3039 & 2.10 & 0.50 & 0.70 & 165.35    & 7.79     & 27         && 3.55     &          & 30        && 2.56         &              & 30             \\
j3040 & 2.10 & 0.50 & 1.00 & 21.78     & 8.30     & 24         && 1.29     &          & 30        && 1.14         &              & 30             \\
j3041 & 2.10 & 0.75 & 0.20 & 1200      & 7.37     & 1          &&          & 26.73    & 0         && 988.73       & 21.38        & 1              \\
j3042 & 2.10 & 0.75 & 0.50 & 258.48    & 11.13    & 10         && 141.39   & 10.73    & 26        && 55.43        & 3.80         & 27             \\
j3043 & 2.10 & 0.75 & 0.70 & 258.89    & 7.82     & 12         && 57.31    & 1.25     & 27        && 30.39        &              & 30             \\
j3044 & 2.10 & 0.75 & 1.00 & 65.36     & 9.21     & 19         && 1.56     &          & 30        && 1.44         &              & 30             \\
j3045 & 2.10 & 1.00 & 0.20 & 660       & 8.68     & 2          &&          & 34.92    & 0         &&              & 31.91        & 0              \\
j3046 & 2.10 & 1.00 & 0.50 &           & 16.45    & 0          && 219.17   & 5.83     & 7         && 517.81       & 5.99         & 16             \\
j3047 & 2.10 & 1.00 & 0.70 & 87.92     & 11.13    & 9          && 114.75   & 2.60     & 26        && 80.91        &              & 30             \\
j3048 & 2.10 & 1.00 & 1.00 & 24.61     & 6.26     & 26         && 1.33     &          & 30        && 1.59         &              & 30             \\
\hline
      &     &      &     &           &          & 767        &&          &          & 1082      &&              &              & 1142   \\
\hline\hline
\end{tabular}
\caption{Comparison of primal method \citep{bruni2018computational}, basic model (\ref{eqn:compact1})-(\ref{eqn:compact11}), and extended model including warm-start and transitivity constraints.}
\label{table:results_comparison2}
\end{table}

\restoregeometry

\section{Conclusion} \label{section:conclusion}

This paper has introduced a new mixed-integer linear programming formulation for the robust counterpart to the two-stage adjustable robust RCPSP. This new compact formulation has been derived by considering a reformulation of the second-stage adversarial sub-problem of maximising the worst-case delayed makespan for a project without resource conflicts. The reformulation of this sub-problem is equivalent to a longest-path problem over an augmented project network made from multiple copies of the original project network. Hence, the dual of this longest-path problem can be inserted into the first-stage resource allocation problem to obtain a compact minimisation problem for the full two-stage robust RCPSP.  

The performance of this new formulation has been examined over 1440 instances of varying characteristics and difficulty, and results show that the proposed formulation can be solved by standard optimisation software significantly faster than the current best algorithm for solving this problem, and can be solved to optimality for almost 50\% more instances.

Regarding future research on the two-stage robust RCPSP, the development of heuristic approaches for solving larger and more-challenging instances of this problem would seem to be a natural and worthwhile objective.

\section*{Acknowledgements}

The authors are grateful for the support of the EPSRC-funded (EP/L015692/1) STOR-i Centre for Doctoral Training. 

\singlespacing

\bibliography{paper}

\onehalfspacing

\appendix

\section{Non-integrality of the adversarial sub-problem}\label{sec:appendix}

Here, we show that the constraint matrix of model (\ref{eqn:linearised1})-(\ref{eqn:linearised10}) is not totally unimodular, contrary to the claim made in \cite{bruni2017adjustable}. In the following, we define $\mathcal{E}:=E\cup X$. The constraint matrix of (\ref{eqn:linearised1})-(\ref{eqn:linearised10}) can be written in matrix notation as:

\[
C=
\begin{blockarray}{cccl}
\alpha & w & \delta &\\
\begin{block}{(ccc)l}
	A & 0 & 0 & \textnormal{Group 1 (\ref{eqn:linearised2})-(\ref{eqn:linearised4})}\\
	0 & I_\mathcal{E} & -B & \textnormal{Group 2 (\ref{eqn:linearised5})}\\
	-I_\mathcal{E} & I_\mathcal{E} & 0 & \textnormal{Group 3 (\ref{eqn:linearised6})}\\
	0 & 0 & e^T_V & \textnormal{Group 4 (\ref{eqn:linearised7})}\\
	0 & 0 & I_V & \textnormal{Group 5 (\ref{eqn:linearised8})}\\
\end{block}
\end{blockarray}
 \]
 where $A$ is a $|V|\times|\mathcal{E}|$ arc-node incidence matrix, $B$ is a $|\mathcal{E}|\times|V|$ matrix where $B_{(i,j),i}=1$ for each $(i,j)\in \mathcal{E}$, and 0 otherwise, $I_V$ and $I_\mathcal{E}$ are identity matrices of dimension $|V|$ and $|\mathcal{E}|$ respectively, and $e^T_V$ is a $|V|\times 1$ vectors of 1's. The rows of this matrix have been grouped according to the constraints that they represent, and similarly, the columns have been grouped by the variables that they represent.

\cite{ghouila1962caracterisation} showed that a matrix $A$ is totally unimodular if and only if for every subset of rows $R$, there exists a partition of $R$ into two disjoint subsets $R_1$ and $R_2$ such that 
 $$\sum_{i\in R_1} a_{ij}-\sum_{i \in R_2}a_{ij}\in\{-1,0,1\},\quad \forall j=1,\dots,n.$$
 Therefore, finding a subset of rows of matrix $C$ for which this condition cannot hold will prove that $C$ is not totally unimodular. 
 
Consider the constraint matrix of the example shown in Figure \ref{fig:counter_example}:

\begin{center} 
\includegraphics[scale=0.9]{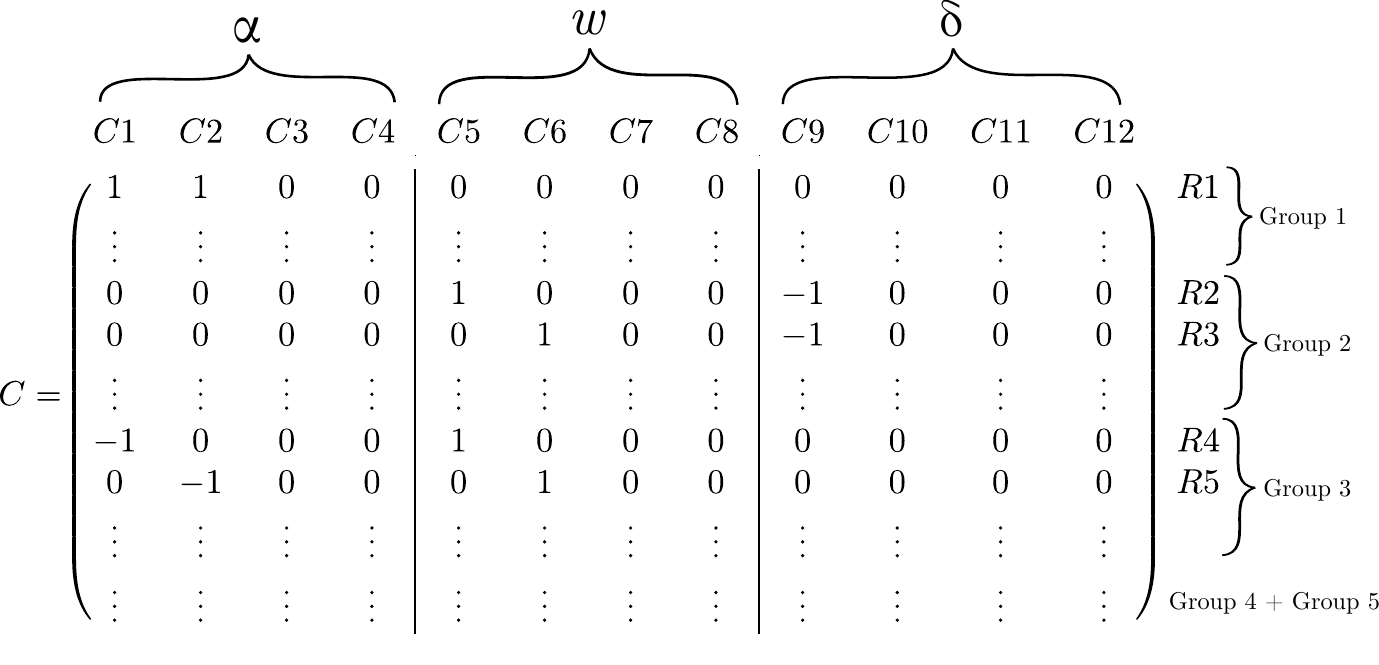}
\end{center}
Take $R$ to be the subset of rows consisting of the first row of Group 1, and first two rows of Groups 2 and 3. We will refer to these rows as $R1,\dots,R5$. To ensure that the sum of column $C9$ is in $\{-1,0,1\}$, $R2$ and $R3$ must be assigned opposite signs. $R4$ and $R5$ must have opposite signs to $R2$ and $R3$, respectively to ensure that the sum of columns $C5$ and $C6$ are in $\{-1,0,1\}$. Then, whatever the choice of sign for $R1$, the sum of column $C1$ and the sum of column $C2$ cannot both be in $\{-1,0,1\}$. Hence, there exists a subset of rows for which the Ghouila-Houri characterisation of total unimodularity does not hold, thus proving that matrix $C$ is not totally unimodular, and that model (\ref{eqn:linearised1})-(\ref{eqn:linearised10}) is not equivalent to its linear relaxation.

\end{document}